\documentclass[%
 reprint,
 amsmath,amssymb,
 aps,
]{revtex4-2}

\usepackage{graphicx}% Include figure files
\usepackage{dcolumn}% Align table columns on decimal point
\usepackage[utf8]{inputenc}
\usepackage[T1]{fontenc}
 \usepackage{xcolor}

\makeatletter
\def\@email#1#2{%
 \endgroup
 \patchcmd{\titleblock@produce}
  {\frontmatter@RRAPformat}
  {\frontmatter@RRAPformat{\produce@RRAP{*#1\href{mailto:#2}{#2}}}\frontmatter@RRAPformat}
  {}{}
}%
\makeatother

\usepackage{amsmath}

\begin{document}

\title{Singular basins in multiscale systems: tunneling between stable states}
% Force line breaks with \\
\author{S. Yanchuk$^{1, 2, *}$}
\author{S. Wieczorek$^1$}%
\author{H. Jard\'on-Kojakhmetov$^3$}
\author{H. Alkhayuon$^{1,}$}
\thanks{These authors contributed equally to this work.}
\affiliation{$^1$School of Mathematical Sciences, University College Cork, Western Road, Cork T12 XF62, Ireland}
\affiliation{$^2$Potsdam Institute for Climate Impact Research, P.O Box 6012 03, 14412 Potsdam, Germany}
\affiliation{$^3$Johann Bernoulli Institute for Mathematics and Computer Science,
University of Groningen, P.O. Box 407, 9700 AK, Groningen, The Netherlands}

\date{\today}% It is always \today, today,

\begin{abstract}
Real-world systems often evolve on  different timescales and possess  multiple coexisting stable states. Whether or not a system returns to a given stable state after being perturbed away from it depends on the shape and extent of its basin of attraction.
We show that basins of attraction in multiscale systems can exhibit  special geometric properties in the form of  \textit{singular funnels}. Although singular funnels are narrow, they can extend to different regions of the phase space and, unexpectedly, impact the system's resilience to perturbations. Consequently, singular funnels may prevent common dimensionality reductions in the limit of large timescale separation, such as the quasi-static approximation, adiabatic elimination and time-averaging of the fast variables.
We refer to basins of attraction with singular funnels as \textit{singular basins}. We show that singular basins are universal and occur robustly in a range of multiscale systems: the normal form of a pitchfork bifurcation with a slowly adapting parameter, an adaptive active rotator, and an adaptive network of phase rotators. 
\end{abstract}

\maketitle

Complex real-world systems are characterized by \textit{multistability}, i.e., when started from different initial conditions, the same system can  end up with notably different asymptotic behavior.
Examples of multistability can be found in epileptic and neuronal models \cite{Gerstner2014,Pfeifer2021},
reservoir computers \cite{Flynn2021a}, lasers 
\cite{Jaurigue2017,Wieczorek2005}, and climate \cite{alkhayuon2019basin, slyman2025tipping}, to name just a few. 
The analysis and control of multistability-related phenomena has therefore been the subject of many studies \cite{Pisarchik2014,Malashchenko2011,Hellmann2018}.
An important question in multistable systems concerns possible transitions between different stable states.
One notable approach addressing the above challenges is the concept of basin stability \cite{Menck2013, soliman_global_1992}, which represents the likelihood of reaching a given  state (attractor) when starting from a random initial condition.
This likelihood can be estimated by the relative volume of its basin of attraction, i.e., the proportion of initial states in the phase space that converge to this attractor.
More importantly, the shape of the basins of attraction determines the uncertainty~\cite{Daza2016} and resilience~\cite{holling_resilience_1973, krakovska_resilience_2024} of a system to perturbations. 
Therefore,  to understand and control a multistable system, it is crucial to understand the geometry of its basins of attraction.

In addition to having multiple stable states, it is also common for real-world systems to evolve on {\it multiple timescales}, which gives rise to a variety of nonlinear phenomena \cite{Kuehn2015,Krupa1997,balzer_canard_2024,ritchieRateinducedTippingNatural2023,Desroches2012}.
An important example of a multi-scale system is an adaptive dynamical network, where the adaptation is much slower than the node dynamics \cite{Berner2023,Berner2021,Maslennikov2018,gross2008adaptive,Aoki2011,venegas-pinedaCoevolutionaryControlClass2025}. Such systems are challenging because they are also high-dimensional and can exhibit a high degree of multistability. The study of basins of attraction in multiscale systems is an area of research that remains largely unexplored. Therefore, uncovering new properties of the basins that arise from the interplay of different timescales is an intriguing and important task. 

Useful and widely used techniques have been developed to analyze mathematical models of multiscale systems with two distinct timescales, also known as slow-fast  systems~\cite{jones_geometric_1995,wechselberger_geometric_2020,Fenichel1979}. 
These techniques exploit the limit of infinite ratio between the fast and slow timescales and can be broadly divided into two groups~\cite{hummel_reduction_2023,hasselmann_stochastic_1976,Kuehn2015}. The first group includes adiabatic elimination and averaging. It eliminates the fast variables by treating them as instantaneous, or by taking into account their effective contribution, while only considering the time evolution of the slow variables in the {\em reduced  system}.
The second group includes the quasistatic approximation. It `freezes' the slow variables by treating them as static at each instant, while only considering the time evolution of the fast variables in the {\em layer system}~\cite{Kuehn2015,wechselberger_geometric_2020}.

This Letter uncovers special properties of basins of attraction in slow–fast systems, leading to counter-intuitive phenomena. 
We show that such basins can contain \textit{singular funnels} (SFs) in the form of tunnels that extend to different regions of the phase space and become increasingly narrow as the timescale ratio increases.  
The problem with SFs is that they are eliminated by the quasistatic approximation and adiabatic elimination. Consequently, adiabatic elimination makes it impossible to reach a given stable state from certain regions of the phase space, whereas this is possible in the full system.
On the other hand, the quasistatic approximation could mask the presence of multiple basins of attraction entirely. 
Therefore, one must be cautious when extrapolating the resilience properties of the full system from those of a reduced one, regardless of whether the latter is obtained via quasistatic approximation, adiabatic elimination, or averaging.
Due to the vanishing of SFs in the limit of infinite ratio of the
timescales, we refer to basins containing SFs as \textit{singular basins}.

\begin{figure}
    \centering
    \includegraphics[width=1\linewidth]{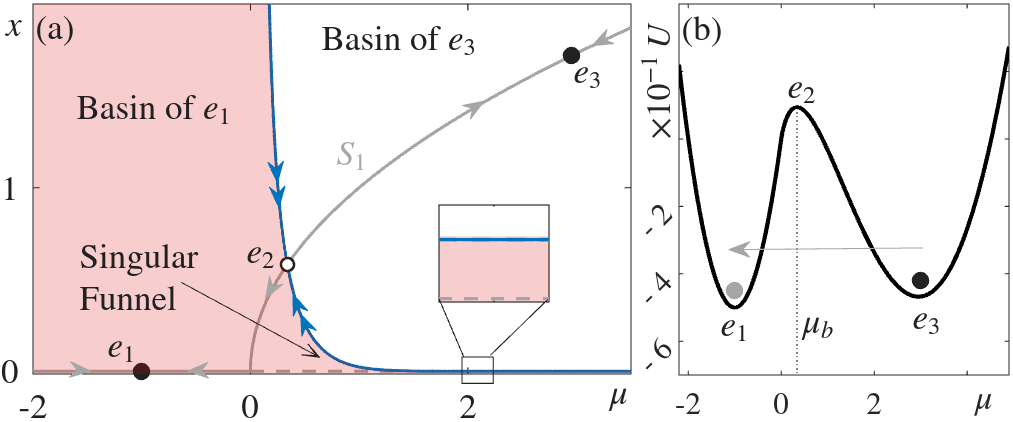}
    \caption{
    % {\sy (HA TODO)}
    % {\sw\em (The dashed part of the grey line at zero is not visible in (a). Could we start the $x$-axis just below zero or lift the grey curve at zero a bit - either should resolve the problem.)}
    (a) Singular basin with a singular funnel (SF) in a bistable pitchfork normal form (\ref{eq:normal-form-a}--\ref{eq:normal-form-b}) with slowly adapting parameter $\mu$ shown on the horizontal axis. The red (white) region is the basin of attraction of the stable equilibrium $e_1$ ($e_3$). 
    The boundary of the two basins is given by the stable manifold of the saddle equilibrium $e_2$ (blue curves). 
    The SF near $x=0$ shrinks exponentially with $\mu$. 
    (b) The bistable potential of the corresponding slow subsystem after adiabatic elimination of the fast variable $x$ (equivalently, reduction to the stable critical manifold). Parameters:  $a=3$ , $b=2$, and $\varepsilon=0.1$.
     % (c)-(d) The same as (a)-(b) but for system \eqref{eq:new_normal-form-mu} with  $a=5$ , $b=10$, and $\varepsilon=0.1$.
    \label{fig:normal-form}
    }
\end{figure}

We begin by illustrating the surprising effects of SFs in the normal form for a supercritical pitchfork bifurcation with a slowly adapting bifurcation parameter. The normal form alone is given by $\frac{dx}{dt} = x(\mu-x^2)$ \cite{Kuznetsov1995}. We restrict to $x \ge 0$, which is appropriate for many applications, for example, when the variable $x$ denotes a population density or the amplitude of a certain observable. Such system has a single stable equilibrium, $x^*$, which exists continuously for all values of $\mu$. Specifically, $x^*=0$ for $\mu\le0$, and $x^*=\sqrt{\mu}$ for $\mu > 0$, 
so it can be written in a more compact form as
$x^*(\mu) =  \sqrt{\mu H(\mu)}$, where $H(\mu)$ is the Heaviside step function \footnote{$H(\mu)=0$ for $\mu<0$ and $H(\mu)=1$ for $\mu\ge 0$.}.
We then introduce the slow variable in the form of a linear adaptation of the bifurcation parameter $\mu$, resulting in the following slow-fast system
written in terms of the {\em fast time} $t$,
\begin{align}
        \frac{\mathrm{d} x}{\mathrm{d}t} & = x(\mu-x^2),     \label{eq:normal-form-a} \\
        \frac{\mathrm{d} \mu}{\mathrm{d}t} &= \varepsilon (- \mu +a x - b). \label{eq:normal-form-b}
\end{align}
The small parameter $0 <\varepsilon \ll 1$ is the ratio between 
the slow timescale of $\mu$ and the fast timescale of $x$.
When $\varepsilon = 0$, the slow timescale is infinitely slower, or equivalently, the fast timescale is infinitely faster. The parameters $a,b>0$ determine the adaptation rule. Under the condition $a>2\sqrt{b}$, system (\ref{eq:normal-form-a}--\ref{eq:normal-form-b}) has two stable equilibria, $e_1$ and $e_3$, one saddle equilibrium, $e_2$, and an SF, as shown in Fig.~\ref{fig:normal-form}(a); see \cite{supplem} for more details.

The layer system is obtained by setting $\varepsilon=0$ in (\ref{eq:normal-form-a}--\ref{eq:normal-form-b}), which gives system \eqref{eq:normal-form-a} with  a static $\mu$. 
It has a unique stable quasistatic equilibrium, $x^*(\mu)$, which varies continuously with $\mu$ and attracts all initial conditions $x_0>0$ at every value of $\mu$ (the solid part of the gray curve $S$ in Fig.~\ref{fig:normal-form}(a)).
Therefore, we apply the adiabatic elimination procedure.
% This suggests that it should be safe to 
% apply the adiabatic elimination procedure globally.
We therefore substitute  $x^*(\mu)$ into \eqref{eq:normal-form-b}, define the {\em slow time} to be $\tau =\varepsilon t$, and obtain the reduced system in terms of $\tau$,
\begin{equation}
    \label{eq:nf-reduced}
    \frac{\mathrm{d}\mu}{\mathrm{d}\tau} = f(\mu),
\end{equation}
where $f(\mu)=-\mu +a  \sqrt{\mu H(\mu)} - b$. 
%and $\tau =\varepsilon t$ {\sw is} the slow time. 
We then recognise that this is a bistable dynamical system with a double-well potential, $U(\mu) = -\int f(\mu)d\mu = \mu^2/2 +b\mu -2 a \left(\mu H(\mu)\right)^{3/2} /3$, and the same equilibria as the full system; 
see Fig.~\ref{fig:normal-form}(b).
After adiabatic elimination, the basins of attraction of the stable equilibria $e_1$ and $e_3$ in Fig.~\ref{fig:normal-form}(b) are clearly separated by a single value of the slow variable, $\mu =\mu_b$, given by the position of the unstable equilibrium $e_2$. If the reduced
system 
is  started on the left (resp. right) of $\mu_b$, it can only converge to $e_1$ (resp. $e_3$). 
However, there is a surprise. This clear separation of basins by $\mu=\mu_b$ is never present in the full system 
(\ref{eq:normal-form-a}--\ref{eq:normal-form-b}) due to the SF, no matter how large the timescale ratio is (or how small but non-zero the value of $\varepsilon$ is). In particular, the SF allows convergence to $e_1$ from both sides of $\mu_b$, which is impossible in the reduced system.
This includes perturbations from $e_3$ into the SF along the fast $x$-direction.  From the  reduced system's perspective, this appears as though there is quantum-like `tunnelling' or `teleporting' from $e_3$ to $e_1$ (the grey arrow in Fig.~\ref{fig:normal-form}(b)). 
In summary, the SF clearly restricts the applicability of the adiabatic elimination,  even though
the layer system has a unique stable quasistatic equilibrium
which varies continuously with $\mu$ and attracts all initial conditions $x_0>0$ at every value of $\mu$.
To demonstrate
that the failure of a global adiabatic elimination for system (\ref{eq:normal-form-a}--\ref{eq:normal-form-b}) is not caused by the degeneracy at 
% {\sy (I would not add it)} {\sw the pitchfork bifurcation point} 
$(\mu,x)=(0,0)$, the supplemental material \cite{supplem} presents a modified system that does not have such degeneracy, but still has an SF.

The SF in the paradigmatic model (\ref{eq:normal-form-a}--\ref{eq:normal-form-b}) has a very simple geometry. In general, however, SFs can have intricate geometries and can penetrate the (white) basin of the other stable state in complicated ways. We will now demonstrate this for adaptive phase rotators. 

The fast-slow adaptive phase rotator has the form \cite{Franovic2020} 
\begin{align} 
        \frac{\mathrm{d}{\varphi}}{\mathrm{d}t}&=\omega+\mu-\sin\varphi, \label{eq:sing_osc-a}\\
        \frac{\mathrm{d}{\mu}}{\mathrm{d}t}&=\varepsilon(-\mu +\eta(1-\sin(\varphi+\alpha))),\label{eq:sing_osc-b}
\end{align}
where $\varphi \in (0,2\pi]$ is a phase variable, $\omega$ is the base oscillator frequency, and $\mu$ accounts for a slow self-adjustment of the oscillator frequency. 
The layer system \eqref{eq:sing_osc-a} with static $\mu$ has a pair of  equilibria, one attracting and one repelling, for $|\omega+\mu|<1$, and a periodic rotation for $|\omega+\mu|>1$ \footnote{This corresponds to the case of unboundedly growing $\varphi(t)$ when $\varphi$ considered on a real line $\mathbb{R}$ or, equivalently, rotating on a circle $\mathbb{S}^1$ when $\varphi$ considered modulo $2\pi$ phase}. At the point $|\omega+\mu|=1$, the layer  system undergoes a saddle-node on invariant circle bifurcation \cite{Strogatz1994, Izhikevich2007,Ermentrout2010} \footnote{At this bifurcation, 
%a saddle and a node 
two equilibria, one stable and one unstable, collide on a closed invariant curve that rotates around a cylinder  due to periodic phase variable. Past the bifurcation, the closed invariant curve becomes a limit cycle whose period tends to infinity at the bifurcation point.
%it becomes a rotating periodic orbit, with period tending to infinity {\sw as the bifurcation point is approached.}
%at onset.
This bifurcation has been referred to as the 'saddle node on invariant circle' (SNIC) bifurcation in \cite{Izhikevich2007}, the 'saddle-node on a limit cycle' (SNLC) bifurcation in \cite{Ermentrout2010}, and the 'infinite-period bifurcation' in \cite{Strogatz1994}}. 
Hence, the layer system \eqref{eq:sing_osc-a} has a unique attracting quasistatic equilibrium, which varies continuously with $\mu$ and attracts all initial conditions (except for the repelling equilibrium point) at every value of $\mu$, when $|\omega+\mu|\le 1$. The attracting equilibrium is replaced by a periodic rotation that exists when $|\omega+\mu|> 1$.

\begin{figure}
    \centering
    \includegraphics[width=1\linewidth]{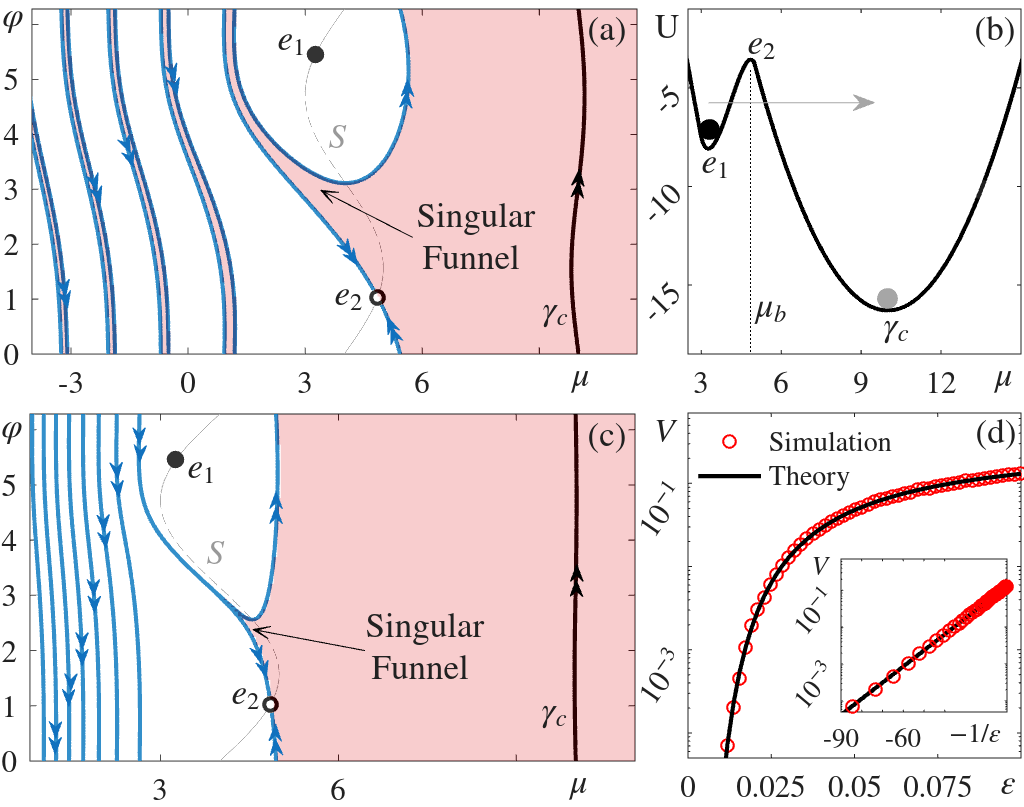}
    \caption{ 
    (a,c) Singular basin with a singular funnel (SF) in an adaptive phase-oscillator (\ref{eq:sing_osc-a}--\ref{eq:sing_osc-b}). The basin of attraction of the periodic rotation $\gamma_c$ is plotted in red and the basin of attraction of the stable equilibrium $e_1$ is plotted in white.
    The (blue) boundary of the two basins is given by the stable manifold of the saddle equilibrium $e_2$. 
    (b) The bistable  potential of the corresponding slow subsystem after adiabatic elimination and averaging of the fast variable $\varphi$.
    (d) The volume of the SF limited to $-10 \le \mu\le 3$ as a function of $\varepsilon$ obtained using equation \eqref{eq:scaling} (black curve) and Monte Carlo simulation (red open circles) \cite{supplem}. 
    Other parameters: $\eta = 10$, $\omega = -4$; (a): $\alpha = \pi/2$, $\varepsilon = 0.1$; (c): $\alpha = \pi/2$, $\varepsilon = 0.01$.
    \label{fig:oneosc} 
    } 
\end{figure}

The adaptive phase rotator (\ref{eq:sing_osc-a}--\ref{eq:sing_osc-b}) can also be reduced to the form \eqref{eq:nf-reduced} by eliminating the fast variable $\varphi$. This is achieved by averaging over the fast rotations for $|\omega+\mu|>1$, and by adiabatic elimination for $|\omega+\mu|<1$.
The corresponding function $f(\mu)$ of the reduced system \eqref{eq:nf-reduced} can be calculated explicitly \cite{supplem}:
\begin{equation}
\label{eq:oneosc-reduced}
    f(\mu)=-\mu+\eta\left(1-\left(\mu+\omega\right)\cos\alpha+\Omega(\mu)\right),
\end{equation}
where
\begin{equation*}
    \label{eq:Omega}
    \Omega(\mu)=
    \begin{cases}
        -\sqrt{1-\left(\mu+\omega\right)^{2}}\cdot \sin\alpha \ \ &\text{for} \ \left|\mu+\omega\right|\le 1,\\
        \sqrt{\left(\omega+\mu\right)^{2}-1}\cdot \cos\alpha \ \ &\text{for} \ \mu+\omega > 1,\\
        - \sqrt{\left(\omega+\mu\right)^{2}-1}\cdot \cos\alpha \ \ &\text{for} \ \mu+\omega <-1.
    \end{cases}
\end{equation*}
We consider the case when system (\ref{eq:sing_osc-a}--\ref{eq:sing_osc-b}) possesses two coexisting attractors:  
a stable equilibrium $e_1$ and a stable rotation $\gamma_c$; see Fig.~\ref{fig:oneosc}(a)
and \cite{supplem} for more details.
The reduced system \eqref{eq:nf-reduced}, with $f(\mu)$ given by equation \eqref{eq:oneosc-reduced}, is also a bistable system with a double-well potential as shown in Fig~\ref{fig:oneosc}(b).

The singular basin of the adaptive phase oscillator (\ref{eq:sing_osc-a}--\ref{eq:sing_osc-b})  in Fig.~\ref{fig:oneosc}(a,c) has an SF that
extends into the region of negative $\mu$. This creates a channel through which initial conditions with arbitrary $\mu$ and selected $\varphi$ are attracted to the periodic rotation $\gamma_c$. Once again, this is impossible within the framework of the reduced
system~\eqref{eq:nf-reduced}, where only initial conditions with $\mu > \mu_b$ can converge to the averaged periodic rotation $\gamma_c$; see Fig.~\ref{fig:oneosc}(b).

We emphasize here that the geometry of the SF of the adaptive phase oscillator (\ref{eq:sing_osc-a}--\ref{eq:sing_osc-b}) is more complex than that of the adaptive pitchfork normal form (\ref{eq:normal-form-a}--\ref{eq:normal-form-b}). The main differences are that: (i) The SF extends over the entire range of the phase variable $\varphi$ in Fig.~\ref{fig:oneosc},
as opposed to an exponentially small region of the $x$-variable in Fig.~\ref{fig:normal-form}.
(ii) The SF in 
Fig.~\ref{fig:normal-form} decreases exponentially with $\mu$ with $\mu$ increasing up to infinity, whereas in Fig.~\ref{fig:oneosc}, such an exponential decrease only occurs over a finite $\mu$-interval. This leads to a higher volume of the SF for the system (\ref{eq:sing_osc-a}--\ref{eq:sing_osc-b}) compared to the system (\ref{eq:normal-form-a}--\ref{eq:normal-form-b}).
To provide a better insight into the geometry of the SF for (\ref{eq:sing_osc-a}--\ref{eq:sing_osc-b}), we have included the phase portrait with the variable $\varphi$ shown on the real line in supplemental material \cite{supplem}.

Next, we provide details of how the boundaries of the SF in Fig.~\ref{fig:oneosc}(a,c) are formed. As in the case of the adaptive pitchfork normal form (\ref{eq:normal-form-a}--\ref{eq:normal-form-b}), these boundaries consist of orbits attracted to the saddle equilibrium $e_2$ (branches of the stable  invariant manifold of $e_2$). Tracking these orbits backwards in time reveals that they pass close to the branch of unstable quasistatic equilibria of the layer  system.
This is the dashed part of the curve $S$ in Fig.~\ref{fig:oneosc}(a,c), where $S$ is also known as the critical manifold \cite{Kuehn2015}. Since this part of the critical manifold is attractive in backward time, the two SF boundaries exponentially approach each other as they are attracted to $S$, causing the SF to become exponentially narrow. 
After the SF extends past $S$ to negative $\mu$, it undergoes further rotations,   repeatedly intersecting the white basin of attraction of $e_1$ in
%leading to an SF with 
a complicated pattern. 

For the systems with a singular basin that have been considered so far, we can obtain the universal scaling of 
%the SF volume in phase space as a function of $\varepsilon$ to estimate 
how
%quickly
the SF volume shrink as $\varepsilon$ approaches 0.
In all cases, the motion along the trajectories which correspond to the SF boundaries is slow and its speed is proportional to $\varepsilon$ along the unstable branch of a critical manifold $S$. Thus, the time they spend close to an unstable branch of $S$ of length $L$ scales as $L/\varepsilon$. The width $\delta$ of the SF is determined by two trajectories that pass close to $S$ (blue trajectories in Fig.~\ref{fig:oneosc}). Hence, $\delta$ decreases exponentially as time goes backwards. Therefore,  this width can be estimated as $\delta \sim \exp(-L\lambda /\varepsilon)$ at the moment when the SF leaves the neighbourhood of an unstable branch of $S$ of length $L$. 
Here $\lambda$  is an effective repulsion rate  away from
%at 
$S$, which corresponds to the attraction in reverse time. 
These arguments lead to the following expected scaling for the volume of an SF:  
\begin{equation}
\label{eq:scaling}
V(\varepsilon) \sim \exp(- C \varepsilon^{-1}), \quad C>0. 
\end{equation}
The scaling \eqref{eq:scaling} holds for the adaptive normal form (\ref{eq:normal-form-a}--\ref{eq:normal-form-b}), since the stable manifold of $e_2$ is exponentially repelled from $x=0$. Furthermore, the adaptive phase rotator (\ref{eq:sing_osc-a}--\ref{eq:sing_osc-b}) exhibits this scaling too, as shown in Fig.~\ref{fig:oneosc}(c). 
The SF scaling \eqref{eq:scaling} is expected to hold for large classes of slow-fast systems, at least in low dimensions, as the arguments leading to \eqref{eq:scaling} are based on the rather general geometric structure of the basin boundaries. 
However, as we show below,  SFs with  even more intricate geometries found in higher-dimensional systems may show deviations from this scaling.

\begin{figure}[t]
    \centering
    \includegraphics[width=1\linewidth]{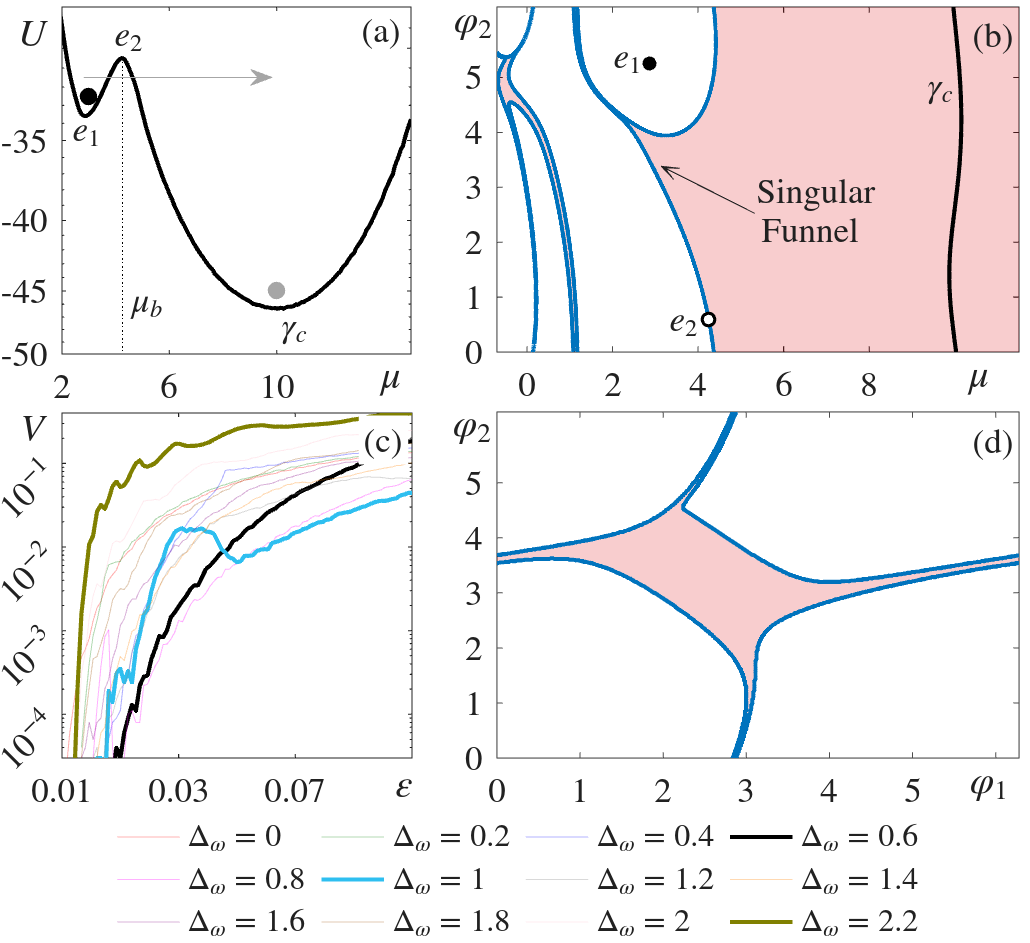}
    \caption{
    (a) Bistable potential of the slow subsystem of two adaptive phase oscillators (\ref{eq:N_osc}-\ref{eq:mudot}), for $N = 2$, after adiabatic elimination and averaging of the fast variables $\varphi_{1,2}$.    
    (b, d) Two cross-sections of the singular basin with a singular funnel (SF) of the same system: (b) $\varphi_1 = 1.2461$ and (d) $\mu = 2.86$. 
    The basin of attraction of the periodic rotation $\gamma_c$ is plotted in red and the basin of attraction of the stable equilibrium $e_1$ is plotted in white.
    The (blue) boundary of the two basins is given by the stable manifold of the saddle equilibrium $e_2$. 
    (c) 
    The volume of the SF limited to $-10 \le \mu\le 0$ as a function of $\varepsilon$ obtained using Monte Carlo simulation \cite{supplem} for different values of the parameter $\Delta_\omega=\omega_1-\omega_2$, showing the robustness of the SF.     
    Other parameters: $\omega_1 = -4$, $\omega_2 = -3$, $\kappa = 1$, $\eta = 10$, $\alpha = \pi/2$; (b-d) $\varepsilon = 0.1$.
    %\sy{plot legends in a more compact way}
    % Illustration of the singular basin in system \eqref{eq:N_osc} of $N=2$ coupled phase rotators {\sw\em (System (9) does not have a singular basin! It has a regular basin.)}. (a) The bistable potential of the averaged system. (b) A two-dimensional cross-section of the basin of the full slow-fast system. The notations are the same as in Fig.~\ref{fig:oneosc}. The cross-section for fixed $\varphi_1 = 1.2461$.
    % (c) The volume of the singular funnel for $\gamma_c$, limited to $\mu\in [-10,\,0]$, as a function of $\varepsilon$ for different values of the parameter $\Delta_\omega=\omega_1-\omega_2$. (d) The cross-section of the basin for fixed $\mu=2.86$.
    }
    \label{fig:network-2}
\end{figure}

A higher-dimensional class of systems, in which singular basins can be observed, is the system of mean-field coupled active rotators
\begin{align}
    \label{eq:N_osc}
    \frac{\mathrm{d}\varphi_i}{\mathrm{d}t}&=\omega_i + \mu - \sin\varphi_i + \frac{\kappa}{N}\sum_{j=1}^N\sin(\phi_j-\phi_i), \\
    \label{eq:mudot}
    \frac{\mathrm{d}\mu}{\mathrm{d}t}&=\varepsilon(-\mu +\eta(1-X)), 
\end{align}
where $i=1,\dots,N$ and $X =  \frac{1}{N} \sum_{j=1}^N \sin\left(\varphi_j+\alpha\right)$.
In this system, the frequencies of the individual rotators are adapted globally by the slow variable $\mu$, and $\mu$ is driven by the mean-field $X$. 
In contrast to the extensively studied system of coupled rotators without adaptation \cite{Park1996,Lindner2004,Dolmatova2017,Klinshov2021,Ronge2021,burylkoTimereversibleDynamicsSystem2023,Hellmann2018,Schafer2018}, the adaptive system (\ref{eq:N_osc}--\ref{eq:mudot}) evolves on two distinct timescales with ratio $\varepsilon$. Adaptive systems of a similar nature have been shown to exhibit distinct dynamical properties such as canard cascading  \cite{DHuys2021,balzer_canard_2024}, emergent excitability \cite{Ciszak2020a}, and others \cite{Song2020,ciszak_collective_2021,venegas2023stable}. Here, we add singular basins to this list of distinct properties.

Figure \ref{fig:network-2} shows a singular basin and its scaling in a system of $N=2$ rotators. The system features coexistence of two attractors: one stable equilibrium and one periodic or quasi-periodic rotation. 
A bistable potential calculated by averaging system \eqref{eq:N_osc} with $N=2$ oscillators is shown in Figure~\ref{fig:network-2}(a).  Two-dimensional cross-sections of the corresponding singular basin for fixed values of $\varphi_1$ and $\mu$ values are illustrated in Figs.~\ref{fig:network-2}(b) and (d), respectively.  

The following observations shed light on the topological structure of the SF in the three-dimensional model  (\ref{eq:N_osc}--\ref{eq:mudot}) with $N=2$: (i) The boundary of the SF is formed by the stable invariant manifold of the saddle point $e_2$. Evidence for this can be seen in Fig.~\ref{fig:network-2}(b), where the boundary of the basins crosses $e_2$ (see also Fig.~S5 in \cite{supplem}). (ii) The SF is organized around a non-trivial unstable critical set, which determines the star-shaped cross-section of the SF shown in Fig.~\ref{fig:network-2}(d), see section V of the supplemental material \cite{supplem} for more details.

Figure~\ref{fig:network-2}(c) shows how the SF volume shrinks with $\varepsilon$ for different values of the detuning parameter $\Delta_\omega=\omega_1-\omega_2$. We observe that $V(\varepsilon)$ follows the scaling \eqref{eq:scaling} for some values of $\Delta_\omega$ (e.g. $\Delta_\omega=0.6$), while for other values of $\Delta_\omega$, it exhibits resonance-like deviations from \eqref{eq:scaling} (e.g. $\Delta_\omega=1$ or $\Delta_\omega=2.2$). Hence, we conclude that higher-dimensional systems can have SFs that deviate from the simple scaling \eqref{eq:scaling} and leave this as a question for future research.

%Consequently, the corresponding averaged system exhibits a bistable potential, as shown in Fig.\ref{fig:network-2}(a). This bistability is caused by the phase-space structures similar to those in the single adaptive rotator. However, a detailed description of these structures is beyond the scope of this letter. 

\begin{figure}[t]
    \centering
    \includegraphics[width=1\linewidth]{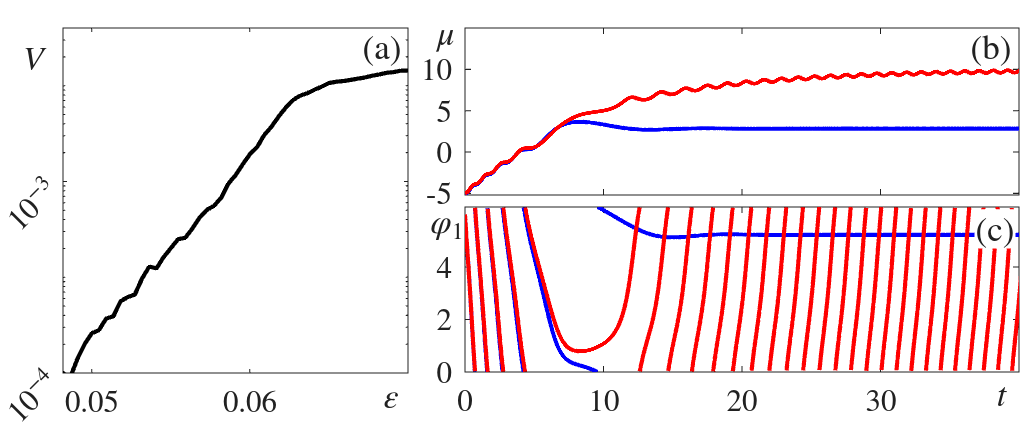}
    \caption{
    (a) The volume of the basin of attraction of ten adaptive phase oscillators (\ref{eq:N_osc}-\ref{eq:mudot}) for $N = 10$, limited to $-10 \le \mu\le 0$, as a function of $\varepsilon$ obtained using Monte Carlo simulation \cite{supplem}. 
    Its decrease with $\varepsilon$ 
     is consistent with
    % indicates 
    the presence of a singular funnel (SF).
    (b,c) Two trajectories starting from close initial conditions:
    $\varphi_i(0) = 6$ for $i = 1,\dots,10$, with $\mu(0) = -5$ (red) and $\mu(0) = -5.1$ (blue).
    The red trajectory starts from SF and, therefore, it is attracted to a stable rotating state.
    Other parameters: 
    $\omega_i = -4 + (i-1)/9$, $i = 1,\dots,10$, $\eta = 10$,
    $\alpha = \pi/2$, and $\kappa = 1$.
%    {\sy{add inset in (a) for the basin}}
    }    
    \label{fig:network-10}
\end{figure}

A visualization of SFs in high-dimensional systems via two-dimensional cross sections is very difficult.
%as one can see in the inset of Figure~\ref{fig:network-10}(a) for the case of $10$ coupled rotators. The SFs can be easily missed in the 11-dimensional phase space. 
Instead, the estimated volume of the SF for the case of $10$ coupled rotators (eleven-dimensional phase space) is plotted in Figure~\ref{fig:network-10}(a). Similar to the lower-dimensional examples, the volume shrinks with $\varepsilon$, indicating the presence of an SF. Additionally, Fig.~\ref{fig:network-10}(b) shows two trajectories: (red) one starting from the SF and converging to periodic rotations, and (blue) one starting and staying very close by, but eventually converging to the stable equilibrium $e_1$, which is another indication of an SF.

 Finally, singular basins are a robust phenomenon. This is due to the fact that its ingredients persist under small parameter changes. This applies to the stable invariant manifolds of 
 %the 
 saddle equilibria and the coexisting attractors. We illustrate the wide parameter region in which the singular basin persists and the possible mechanisms  by which it is destroyed in
 %of its destruction for 
 system (\ref{eq:sing_osc-a}--\ref{eq:sing_osc-b}) in \cite{supplem}.

In summary, we have uncovered singular basins with singular funnels (SFs) in bistable slow-fast systems, and showed that they are robust and ubiquitous. The SFs prevent  useful and widely used dimensionality reductions, even when 
there is a unique quasistatic attractor, which changes continuously with the slow variable(s) and attracts almost all initial conditions at every `frozen' setting of the slow variable.
Although adiabatic elimination or averaging of the fast variables retains the coexistence of different stable states and different basins of attraction in the reduced system, the fast dynamics of the full system can lead to transitions between different stable states via SFs that are not possible in the reduced system.  From the reduced system's perspective, these unexpected dynamics resemble quantum-like `tunnelling' or `teleporting' from one attractor to another. 
We have constructed low-dimensional canonical systems that can be considered  normal forms for a singular basin and showed a universal scaling of SFs with increasing timescale ratio in these systems. Additionally, we have demonstrated that SFs can form intricate structures in the phase space of higher-dimensional systems with slowly adaptive mean-field coupling. A challenge for future research is to understand singular basins in systems with more than two timescales and more than two coexisting stable states. 

Perhaps most intriguingly of all, any mathematical model of a real-world system is necessarily an approximation that omits some fast variables. Thus, our results also provide new insight into the potential for unexpected critical transitions~\cite{Ashwin2012}. For example, transitions between different stable states, which are not captured by the models due to omitting some fast variables, could be triggered in real systems.

\begin{acknowledgments}
The work od S.Y and S.W was supported by Taighde Éireann—Research Ireland (Grant No. FFP-A/12066).  
\end{acknowledgments}

~

{
{\bf Code availability:} The code used to perform the simulations and produce all figures is available in the GitHub repository: \url{https://github.com/hassanalkhayuon/Singular_Funnels}
}

\bibliography{bibliography}

\onecolumngrid

\newpage
\begin{center}
{\huge Supplemental Material}
\end{center}

\section{Details of the numerical methods}

\subsection{The singular funnel volume using Monte Carlo simulations in Fig.~2(d)~and~3(c),\label{sec:numarical_basin}}
We estimated the volume of the singular funnel in Fig.~2~and~3 of the main paper by Monte Carlo simulations as follows:
\begin{itemize}
    \item 
    For each value of $\varepsilon$, we randomly choose a large number $M$ of uniformly distributed initial conditions, with $\varphi_i(0) \in [0,\, 2\pi)$ and $\mu(0) \in [-10,\, 0]$. 
    \item 
    For each of these initial conditions, we numerically solve the initial value problem consisting of system~\eqref{eq:sing_osc-a}--\eqref{eq:sing_osc-b} or \eqref{eq:N_osc}--\eqref{eq:mudot}, in the main paper, and the initial condition, for a long enough time.  
    For the numerical integration we use the MATLAB function \texttt{ode45}, with $\texttt{RelTol} = 10^{-10}$ and $\texttt{AbsTol} = 10^{-10}$\footnote{We provide the \texttt{RelTol} and \texttt{AbsTol} values for reproducibility purposes. Larger tolerances may be used to speed up the calculations, while producing results that are practically indistinguishable.}. The time span for our integration is $[0,\, 10/\varepsilon]$. 
    \item 
    We then determine whether the solution trajectory converges to the stable equilibrium $e_1$ or to the rotating periodic orbit $\gamma_c$. 
    To detect convergence, we use a user-defined event function to stop the integration at time $t_{\mathrm{end}}$ when $\mu(t_{\mathrm{end}}) > 9$. 
    If this condition applies, then the solution has converged to $\gamma_c$ otherwise, it has converged to $e_1$. 
    \item 
    The volume $V$ of the basin of attraction of the limit cycle $\gamma_c$ is then given by
    $$
    V = \frac{K}{M},
    $$
    where $K$ is the number of initial conditions whose solution trajectories converge to $\gamma_c$ over time. 
\end{itemize}

{\bf Remark:} In Fig.~2(a) and~(c), the boundary between the two basins of attraction is given by the one-dimensional stable manifold of the saddle equilibrium $e_2$. We computed this manifold numerically (using \texttt{ode45}) by selecting two initial conditions along the stable eigendirection of $e_2$ and integrating backward in time. 
The integration parameters are the same as those used above.
The same applies for Fig.~1(a), Fig.~\ref{figCM:normal-form}(a), and Fig.~\ref{fig:normal-form-tanh}(a). 

\subsection{Computing the basin of $\gamma_c$ in Figs.~3(b)~and~(d)}

We computed the basin of attraction of the rotating periodic solution $\gamma_c$ for the two adaptive rotator example in two 2D projections: in the ($\mu$, $\varphi_2$)-plane with $\varphi_1 = 1.2461$ (Fig.~3(b)), and in the ($\varphi_1$, $\varphi_2$)-plane with $\mu = 2.86$ (Fig.~3(d)). For these calculations, we use a grid of initial conditions in each plane, as follows:

\begin{itemize}
    \item 
    We fix the parameters $\omega_1 = -4, \, \omega_2 = -3, \, \kappa = 1, \, \eta = 10, \, \alpha = \pi/2$, and $\varepsilon = 0.1$, and choose an equally spaced grid of 100,000 initial conditions for $\big(\mu(0), \varphi_2(0)\big) \in [-1,\, 10] \times [0, 2\pi) $ for Fig.~3(b) and $\big(\varphi_1(0), \varphi_2(0)\big) \in [0,\,2\pi) \times [0, 2\pi)$. 
    \item 
    For each of these initial conditions, we numerically solve the initial value problem consisting of \eqref{eq:N_osc}--\eqref{eq:mudot}, from the main paper, and the initial condition, for a long enough time.  
    For the numerical integration we use the MATLAB function \texttt{ode45}, with $\texttt{RelTol} = 10^{-10}$ and $\texttt{AbsTol} = 10^{-10}$. The time span for our integration is $[0,\, 10/\varepsilon]$. 
    \item 
    We then determine whether the solution trajectory converges to the stable equilibrium $e_1$ or to the rotating periodic orbit $\gamma_c$. 
    To detect convergence, we use a user-defined event function to stop the integration at time $t_{\mathrm{end}}$ when $\mu(t_{\mathrm{end}}) > 8$. 
    If this condition applies, then the solution has converged to $\gamma_c$ otherwise, it has converged to $e_1$. 
\end{itemize}

{\bf Remark:} Similar to Fig.~2(a), the boundary between the two basins of attraction is given by the two-dimensional stable manifold of the saddle equilibrium $e_2$. However, for this example, we do not compute this manifold. The blue boundaries in Fig.~3(b) and (d) are curves separating the two basins in the corresponding cross-sections, plotted using the MATLAB \texttt{contour} function.

\subsection{Numerical averaging in Fig~3(a)}

We computed the potential of the average system of the two oscillator network, Figure~3 of the main paper, as follows: 
\begin{itemize}
    \item 
    For fixed parameter values of $\omega_{1,2}, \, \kappa, \, \eta, \, \alpha$, and a given value of $\mu$, we randomly choose initial values for $\varphi_1(0)$ and $\varphi_2(0)$ from the interval $[0, \, 2\pi)$. 
    \item 
    We consider the (fast) layer problem of system~\eqref{eq:sing_osc-a}--\eqref{eq:sing_osc-b}, in the main paper, which is given as:
    \begin{equation}
        \begin{split}\label{eq:layer_2Osc}
            \frac{d\varphi_i}{dt}=\omega_i + \mu - \sin\varphi_i + \frac{\kappa}{N}\sum_{j=1}^N\sin(\varphi_j-\varphi_i), \quad i = 1,2.
        \end{split}
    \end{equation}
    \item 
    The initial value problem consisting of system~\eqref{eq:layer_2Osc} and the randomly chosen initial condition $\left(\varphi_1(0), \, \varphi_2(0) \right)$ was solved numerically for a sufficiently long time interval $[0,\, 600]$. 
    For the numerical integration, the MATLAB function \texttt{ode45} was used, with $\texttt{RelTol} = 10^{-10}$ and $\texttt{AbsTol} = 10^{-10}$. 
    \item 
    The transient part $t \in [0, 100]$ of the solution was discarded to ensure that the solution has converged to an invariant state. 
    \item 
    We compute $X(t) = \frac{1}{2} \sum_{j =1}^{2} \sin(\varphi_j(t) + \alpha)$ for $t\in [100,600]$ and write the average system as:
    
    \begin{equation}
    \frac{d\mu}{d\tau} = \bar{g}(\mu) = -\mu +  \frac{1}{500}\int_{100}^{600}\eta (1- X(s)) ds. 
    \end{equation}
    
    The definite integral on the right-hand side was computed numerically using the trapezoidal rule  with 5001 mesh grid points. 
    \item
    The potential $U (\mu)$ is given by: 
    $$
    U(\mu) = \int_{0}^{\mu} \bar{g}(s)\, ds.
    $$
\end{itemize}

\section{Supercritical pitchfork normal form with adaptive parameter}
The pitchfork normal form with adaptively changing parameter, as introduced in the manuscript, is given by
\begin{align} 
        \frac{d x}{d t} & = x(\mu-x^2),\quad x\ge 0, \label{eq:normal-form-x}\\
        \frac{d\mu}{d t} &= \varepsilon (- \mu +a x - b), 
        \label{eq:normal-form-mu}
\end{align}
where $\varepsilon>0$ is a small parameter, and $a,b>0$ are parameters determining the linear adaptation function. 

\begin{figure}
    \centering
    \includegraphics[width=0.7\linewidth]{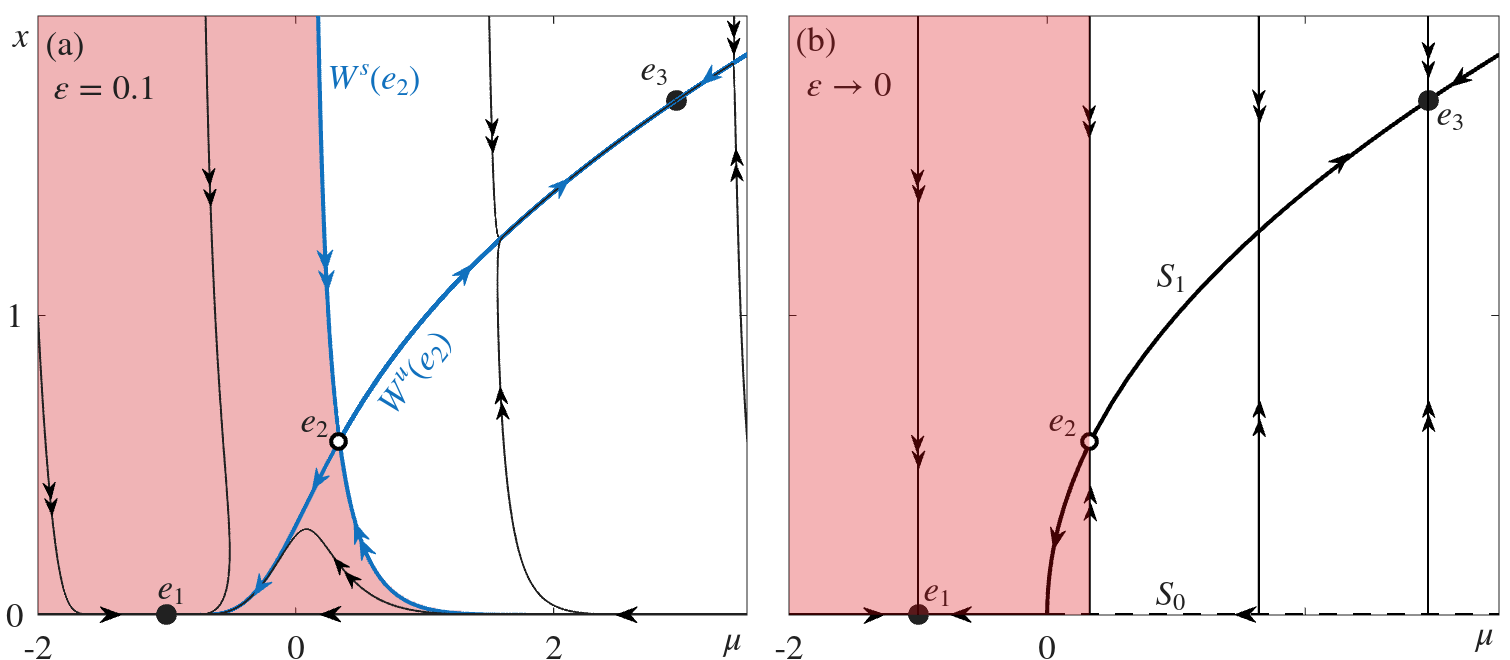}
    \caption{Phase portrait of the pitchfork normal form, system~\eqref{eq:normal-form-x}--\eqref{eq:normal-form-mu}.
    (a) $\varepsilon = 0.1$, the red region (white) region is the basin of attraction of $e_0 (e_2)$, the stable $W^s(e_2)$ (unstable $W^u(e_2)$) manifolds of the saddle point $e_1$ in blue, see Sec.~\ref{sec:numarical_basin} for numerical detailed about computing $W^{s/u}(e_2)$.
    (b) $\varepsilon \to 0$, 
    critical manifolds $S_{0,1}$ in thick black, the attracting parts in solid and the repelling part of $S_0$ is dashed.
    The double arrowed vertical black lines indicate the direction of the fast layer system. 
    Parameters: $a = 3, \, b = 2$.
    % critical manifolds, equilibria, and the vector field wit
    % \hjk{should we use colors for $W^u$ and $W^s$ or is it too much?}\hjk{Should we add the fast fiber at $e_2$?}{\ha Do we like it better with $W^u$ and $W^s$ blue? I prefer the one colour. Also, I looked at the fast fiber at $e_2$, it makes the figure a bit crowded.. HA: Fix the colours } \hjk{Agreed, but then it is ``left to the imagination'' what happens to $W^s$ as $\varepsilon\to0$}
    \label{figCM:normal-form}
    }
\end{figure}

First we consider the equilibria of the fast system \eqref{eq:normal-form-x}, which, when considered as a set in the phase space $(x,\mu)$ of the fast-slow system \eqref{eq:normal-form-x}--\eqref{eq:normal-form-mu}, define the critical manifolds
\begin{align}
    S_0 & = \{ (x,\mu)\in[0,\infty)\times\mathbb R\ : \ x=0 \},\\
    S_1 & = \{ (x,\mu)\in[0,\infty)\times\mathbb R\ : \ x=\sqrt{\mu},\ \mu\ge 0 \}.
\end{align}
Linearization of the vector field along these manifolds shows that $S_0$ has two branches, namely
\begin{align}
    S_0^\textup{a} & = \{ (x,\mu)\in[0,\infty)\times\mathbb R\ : \ x=0,\,\mu<0 \},\\
    S_0^\textup{r} & = \{ (x,\mu)\in[0,\infty)\times\mathbb R\ : \ x=0,\,\mu>0 \},
\end{align}
which are locally attracting and repelling respectively and such that $S_0=S_0^\textup{a}\cup\left\{ 0\right\}\cup S_0^\textup{r}$, while $S_1$ is locally attracting.
This stability information is shown in Fig.~\ref{figCM:normal-form}(b) by solid lines for the stable parts and dashed line for the unstable part of the critical manifolds respectively.

Thus, the fast layer system \eqref{eq:normal-form-x} has a unique stable equilibrium for all values of the slow variable $\mu$, which is attracting all values of $x\ge 0$ for $\mu\le 0$ and all $x>0$ for $\mu>0$. 
Therefore, we consider the reduced slow system on the union of the corresponding stable parts of the critical manifolds. We substitute
\begin{equation}
x = \sqrt{\mu H(\mu)} = 
\begin{cases}
0, & \mu \le 0\\    
\sqrt{\mu}, & \mu >0,
\end{cases}    
\end{equation}
into \eqref{eq:normal-form-mu} and obtain 
\begin{equation}
\frac{d\mu}{dt} 
=
\varepsilon  
\begin{cases}
-\mu -b , & \mu\le  0\\    
-\mu + a \sqrt{\mu} - b , & \mu >0,
\end{cases}    
\end{equation}
or in a more compact form 
\begin{equation}
    \label{eqCM:reduced}
    \frac{d\mu}{dt} = \varepsilon \left(
    -\mu + a \sqrt{\mu H(\mu)} - b
    \right),
\end{equation}
where $H(\mu)$ is the Heaviside step function. 

We consider the case when the reduced system \eqref{eqCM:reduced} has three equilibria $e_0$, $e_1$, and $e_2$ as shown in Fig.~\ref{figCM:normal-form}, their $\mu$ coordinates are $\mu_{0}<0$, $\mu_{1}>0$, and $\mu_{2}>0$, respectively. Requiring three equilibria leads to the following conditions on the parameters: 
\begin{equation}
    \label{eqCM:pars}
    b>0 \quad \text{and} \quad a> 2\sqrt{b}.
\end{equation}
Under the conditions \eqref{eqCM:pars}, the reduced system \eqref{eqCM:reduced} is bistable with two stable equilibria $e_{0}$, $e_2$ and one unstable equilibrium $e_1$: 
\begin{equation}
    \label{eqCM:eq}
    e_1: (0,-b), \quad 
    e_{2,3}: \left( \sqrt{\mu}, \frac{a}{2} 
    \pm \sqrt{\left(\frac{a}{2}\right)^2 -b} \right). 
\end{equation}

The potential of system \eqref{eqCM:reduced} can be calculated as the integral of its right-hand side, leading to
\begin{equation}
    \label{eqCM:potential}
    U(\mu) = \frac{\mu^2}{2} + b\mu -\frac{2}{3}a\mu\sqrt{\mu H(\mu)} ,
\end{equation}
where, for simplicity, we dropped the scaling factor $\varepsilon$. 

The equilibria $e_0$ and $e_2$ are stable in the phase space of the full system \eqref{eq:normal-form-x}--\eqref{eq:normal-form-mu}. 
The attraction basins for these equilibria are separated by the stable manifolds $W_{1,2}^s(e_1)$ of the saddle equilibrium $e_1$, see Fig.~\ref{figCM:normal-form} and Fig.~1(a) of the main manuscript. 
For our purposes, it is worth noting that the manifold $W_2^s(e_2)$ converges exponentially to $x=0$ as $\mu\to\infty$. 
As a result, a part of the attraction basin of the equilibrium $e_0$ extends into an exponentially small region for all $\mu>0$, which we call \textit{singular funnel}. 

\section{A model without critical manifold crossing}

The main manuscript shows the failure of a global adiabatic elimination for system (\ref{eq:normal-form-x}--\ref{eq:normal-form-mu}). However, the critical manifold $x^*(\mu)$ of this system  has a degeneracy at $(\mu,x)=(0,0)$, where $x^*(\mu)$ is stable but not hyperbolic. To demonstrate that this degeneracy is not the reason for the adiabatic elimination failure, we present a modified system here
\begin{equation} \label{eq:new_normal-form-mu}
\begin{aligned}
        \frac{dx}{dt} & = x\left( \tanh\mu + 2 - x\right),\\
        \frac{d\mu}{dt} &= \varepsilon (- \mu + ax - b). 
        % b -> a
        % a -> -b
\end{aligned}
\end{equation}

The critical manifolds are given by 
\begin{align}
    S_0 & = \{ (x,\mu)\in[0,\infty)\times\mathbb R\ : \ x=0 \},\\
    S_1 & = \{ (x,\mu)\in[0,\infty)\times\mathbb R\ : \ x=2 + \tanh\mu \}.
\end{align}
Linearization of the vector field along these manifolds shows that $S_0$ is always locally repelling and $S_1$ is locally attracting. Hence, the corresponding reduced system on the stable critical manifold $S_1$ (adiabatic elimination) is given by
\begin{align}
        \frac{d \mu}{d t} &=  \varepsilon(- \mu - b + a\left(\tanh(\mu) + 2)\right), \label{}
\end{align}
and the corresponding (rescaled) potential $V$ is 
\begin{align*}
        U &= -\int \left(a \tanh(\mu) - \mu - b + 2a\right) \, d \mu \\
        & = 
        \mu^2/2 - (2 a - b) \mu - a\ln(\cosh(\mu)).
\end{align*}

\begin{figure}
    \centering
    \includegraphics[width=0.5\linewidth]{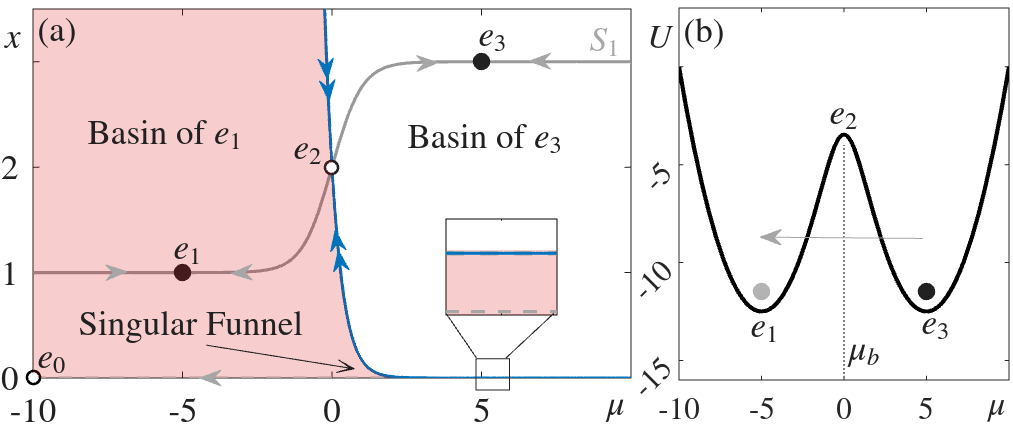}
    \caption{(a) Singular basin with a singular funnel (SF) for system \eqref{eq:new_normal-form-mu}. The red (white) region is the basin of attraction of the stable equilibrium $e_1$ ($e_3$). 
    The boundary of the two basins is given by the stable manifold of the saddle equilibrium $e_2$ (blue curves), see Sec.~\ref{sec:numarical_basin} for numerical details for computing the stable manifold of $e_2$. 
    The SF near $x=0$ shrinks exponentially with $\mu$. 
    (b) The bistable potential of the corresponding slow subsystem after adiabatic elimination of the fast variable $x$ (equivalently, reduction to the stable critical manifold). Parameters:  $a = 5$ , $b = 10$, and $\varepsilon=0.1$.}
    \label{fig:normal-form-tanh}
\end{figure}

In contrast to system \eqref{eq:normal-form-x}--\eqref{eq:normal-form-mu}, the critical manifold of system \eqref{eq:new_normal-form-mu} consists of a single smooth curve $x=2 + \tanh\mu$ without any degeneracies, see Fig.~\ref{fig:normal-form-tanh}. As shown in Fig.~\ref{fig:normal-form-tanh}, this system exhibits a singular basin with an SF and the same failure of global adiabatic elimination.

\section{Adaptive phase rotator}

\subsection{Timescale reduction}

The adaptive phase rotator considered in the manuscript has the form
\begin{align} 
    \frac{d\varphi}{dt}&=\omega+\mu-\sin\varphi, \label{eq:sing_osc-twin-a} \\
    \frac{d\mu}{dt}&=\epsilon(-\mu +\eta(1-\sin(\varphi+\alpha))). \label{eq:sing_osc-twin-b}
\end{align}
The slow-fast dynamics of this system for $\alpha=0$ was studied in \cite{Franovic2020}, for the deterministic and stochastic case. Here we extend these results (for the deterministic case) to $\alpha\ne 0$. 

The critical manifold has the form:
\begin{equation}
    \label{eq:CritMan1osc}
    S_0 := \left\{ (\mu,\varphi):\ \mu=-\omega + \sin\varphi \right\},
\end{equation}
and it exists in the stripe $|\mu+\omega|\le 1$. To obtain the dynamics on this critical manifold (adiabatic elimination), we substitute $\sin\varphi = \mu+\omega $ into \eqref{eq:sing_osc-twin-b}, leading to
\begin{equation}
    \label{eq:SM-adiabatic}
    \frac{d\mu}{dt} =-\mu+\eta\left(1-\left(\mu+\omega\right)\cos\alpha\mp\sqrt{1-\left(\mu+\omega\right)^{2}}\sin\alpha\right),
    \quad 
    \left|\mu+\omega\right|\le1,
\end{equation}
where the negative sign corresponds to $S_0^{\rm a}$ and the positive to $S_0^{\rm r}$.

For $\left|\mu+\omega\right|>1$, the fast subsystem \eqref{eq:sing_osc-twin-a} has no equilibria, but 
 exhibits periodic rotation described by
\[
\varphi_{\mu}(t)=2\arctan\left(\frac{1+\Omega(\mu)\tan\left(\frac{t}{2}\Omega(\mu)\right)}{\omega+\mu}\right),
\]
where 
\[
\Omega(\mu)=\sqrt{\left(\omega+\mu\right)^{2}-1}.
\]

To average the slow dynamics \eqref{eq:sing_osc-twin-b} along these fast rotations, we need to average the oscillating term $\sin(\varphi(t)+\alpha)$ over the period $T=2\pi/\Omega$: 

\[
\left\langle \sin(\varphi_{\mu}(t)+\alpha)\right\rangle =\frac{1}{T}\int_{0}^{2\pi/\Omega}\sin\left(\varphi_{\mu}(t)+\alpha\right)dt=
\]
\[
=\frac{1}{T}\cos\alpha\int_{0}^{T}\sin\varphi_{\mu}(t)dt+\frac{1}{T}\sin\alpha\int_{0}^{T}\cos\varphi_\mu(t) dt=
\]
\[
=
((\omega+\mu)-\textup{sign}(\omega+\mu)\Omega)\cos\alpha
+\frac{1}{T}\sin\alpha\int_{0}^{T}\cos\varphi_{\mu}(t)dt
\]

$$
=
((\omega+\mu)-\textup{sign}(\omega+\mu)\Omega)\cos\alpha
+\frac{1}{T}\sin\alpha\int_{0}^{2\pi}\frac{\cos\varphi d\varphi}{\omega+\mu-\sin\varphi}=
$$
$$
=
((\omega+\mu)-\textup{sign}(\omega+\mu)\Omega)\cos\alpha
+\frac{1}{T}\sin\alpha\int_{0}^{2\pi}\frac{d\sin\varphi}{\omega+\mu-\sin\varphi}=
$$
$$
=((\omega+\mu)-\textup{sign}(\omega+\mu)\Omega)\cos\alpha. 
$$

Therefore, the averaged equation is 
\begin{equation}
\label{eq:SM-averaged}
\frac{d\mu}{dt}=-\mu +
    \eta\left(
    1- 
    \left(
    \left(\omega+\mu\right)-\textup{sign}(\omega+\mu)\sqrt{\left(\omega+\mu\right)^{2}-1}\right)\cos\alpha 
    \right), \quad |\omega+\mu|>1.
\end{equation}

Finally, we combine the reduced systems \eqref{eq:SM-adiabatic} for $|\mu+\omega|\le 1$ and \eqref{eq:SM-averaged} for $|\mu+\omega|> 1$, to obtain
\begin{equation}
    \label{eq:averaged}
    \frac{d\mu}{dt}=-\mu+\eta\left(1-\left(\mu+\omega\right)\cos\alpha+\Omega(\mu)\right),
\end{equation}
where
\begin{equation}
    \label{eq:Omega}
\Omega(\mu)=
    \begin{cases}
        -\sqrt{1-(\omega+\mu)^{2}}\sin\alpha,
        & |(\omega+\mu)|\le 1,\\[1ex]
        \sqrt{(\omega+\mu)^{2}-1}\cos\alpha,
        & (\omega+\mu)>1,\\[1ex]
        -\sqrt{(\omega+\mu)^{2}-1}\cos\alpha,
        & (\omega+\mu)<-1,
    \end{cases}
\end{equation}

The equilibria of the averaged dynamics satisfy
\begin{equation}
    -\mu+\eta\left(1-\left(\mu+\omega\right)\cos\alpha+\Omega(\mu)\right) = 0.
\end{equation}

\subsection{Parameter region for singular basin}

\begin{figure}
    \centering    \includegraphics[width=0.8\linewidth]{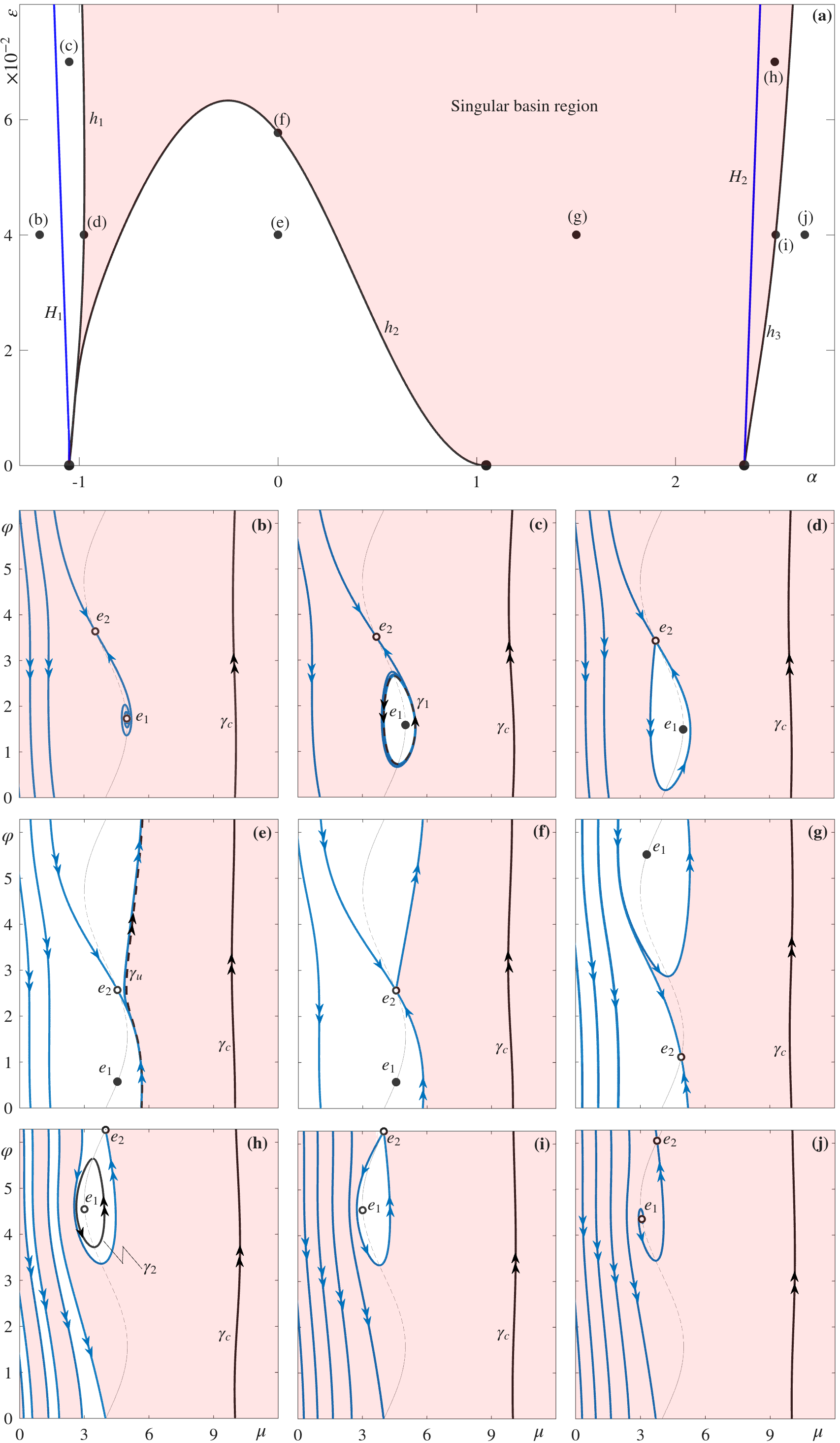}
    \caption{Two-parameter $(\alpha, \, \varepsilon)$ bifurcation diagram for system~\eqref{eq:sing_osc-twin-a}--\eqref{eq:sing_osc-twin-b}, with examples of phase portraits.}
    \label{fig:bif_and_phase_one_osc}
\end{figure}

In Figure~\ref{fig:bif_and_phase_one_osc}, we present the bifurcations of system~\eqref{eq:sing_osc-twin-a}--\eqref{eq:sing_osc-twin-b} with respect to the parameters $\alpha$ and $\varepsilon$.
In the parameter region we examined, there are two Hopf bifurcation curves and three homoclinic curves. Also, the system has two equilibrium solutions $e_1$ and $e_2$, and up to four periodic solutions: $\gamma_{1,2}$, which are regular limit cycles, and $\gamma_{c,u}$, which are rotating limit cycles resulting from the fact that $\varphi \in [0,2\pi)$. We point out that the rotating limit cycle $\gamma_c$ is always stable in the parameter region under examination.

If we consider Figure~\ref{fig:bif_and_phase_one_osc}~(a), starting from the left-hand side,
we have a monostable system in region (b), Figure~\ref{fig:bif_and_phase_one_osc}~(b), where the equilibrium $e_1$ is unstable.
At the subcritical Hopf curve $H_1$, the unstable equilibrium $e_1$ gains stability in region (c), Figure~\ref{fig:bif_and_phase_one_osc}~(c), and an unstable limit cycle $\gamma_1$ emerges to form the boundary of the basin of attraction of $e_1$.

The unstable limit cycle $\gamma_1$ intersects the saddle equilibrium $e_2$ at the homoclinic bifurcation $h_1$, Figure~\ref{fig:bif_and_phase_one_osc}~(d). To the right-hand side of $h_1$, region (g), the basin of attraction of the rotating limit cycle is singular, with the boundary given by the stable manifold of the saddle $e_2$, Figure~\ref{fig:bif_and_phase_one_osc}~(g).

The second homoclinic bifurcation, $h_2$, is formed by the intersection of the unstable rotating limit cycle $\gamma_u$ and the saddle equilibrium $e_2$, Figure~\ref{fig:bif_and_phase_one_osc}~(g).
Region (e), Figure~\ref{fig:bif_and_phase_one_osc}~(e), below the curve $h_2$, is a bistable region where both the equilibrium $e_1$ and the rotating limit cycle $\gamma_c$ are stable, and the basin boundary is given by the unstable rotating limit cycle $\gamma_u$.
The second Hopf bifurcation, $H_2$, is supercritical. The stable equilibrium $e_1$ in Figure~\ref{fig:bif_and_phase_one_osc}~(g) loses stability, giving rise to a stable limit cycle $\gamma_2$, Figure~\ref{fig:bif_and_phase_one_osc}~(h). The basin of attraction of $\gamma_c$ is still singular in region (h).
The region of singular basin ends at the homoclinic bifurcation $h_3$, Figure~\ref{fig:bif_and_phase_one_osc}~(i), where the stable limit cycle $\gamma_2$ intersects the saddle equilibrium $e_2$.
In region (j), Figure~\ref{fig:bif_and_phase_one_osc}~(j), the system is monostable again, where the rotating limit cycle $\gamma_c$ is the only stable attractor.

\begin{figure}
    \centering    \includegraphics[width=0.4\linewidth]{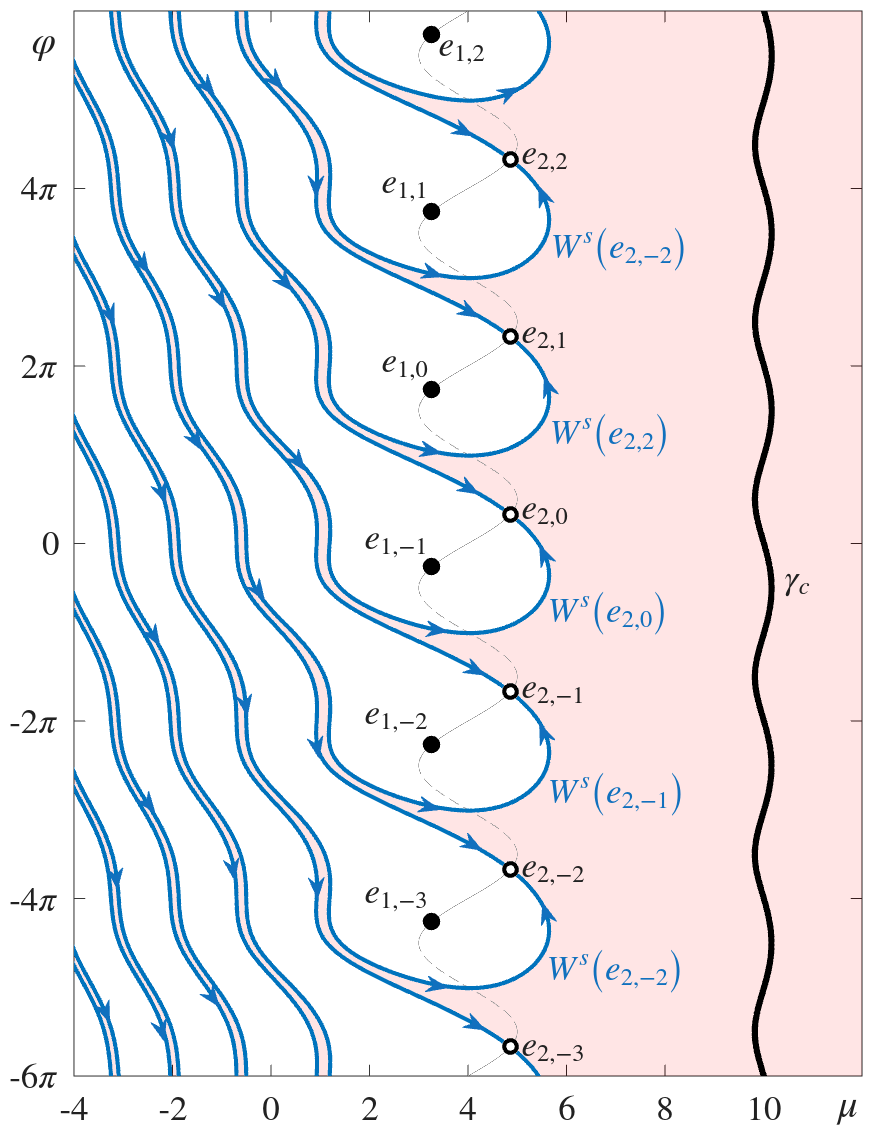}
    \caption{Phase portrait of system~\eqref{eq:sing_osc-twin-a}--\eqref{eq:sing_osc-twin-b} with the phase variable $\varphi \in \mathbb{R}$ unfolded. Parameter values: $\omega = -4, \, \eta = 10, \alpha = \pi/2$ and $\varepsilon = 0.1$. }
    \label{fig:phase_one_osc_unfolded}
\end{figure}

\subsection{Analytical considerations for the singular funnel scaling}

Here we provide some analytical considerations, which substantiate the scaling given in the main manuscript as equation \eqref{eq:scaling}. 

Consider $\varepsilon\ll 1$ and the boundaries of the singular basin to be
\[
\varphi_{u}(t)=\varphi(t,\varphi_{u0},\mu_{0}), \quad 
\mu_{u}(t)=\varphi(t,\varphi_{u0},\mu_{0})
\]
for the upper, and
\[
\varphi_{l}(t)=\varphi(t,\varphi_{l0},\mu_{0}), \quad 
\mu_{l}(t)=\varphi(t,\varphi_{l0},\mu_{0})
\]
for the lower boundaries. Here $\varphi_{l0}<\varphi_{u0}$, $-\omega-1<\mu_0<\mu(e_2)$ are some chosen points on the singular funnel such that $\varphi(0,\varphi_{X0},\mu_0)=\varphi_{X0}$ and $\mu(0,\varphi_{X0},\mu_0)=\mu_{0}$ for $X\in\{u,l\}$. 
For any fixed $\delta<\varphi_{u0}-\varphi_{l0}$, there exists $t_0<0$ such that 
\[
\varphi_{l}(t_{0})-\varphi_{u}(t_{0})=\delta
\]
due to the contraction to the equilibrium of the layer system (fast equation \eqref{eq:sing_osc-twin-a}) in the backwards time. 

Being $\delta$-close to the equilibrium of the fast system, it can be linearized for $t_f<t<t_0$, where $t_f$ is the time when the system approaches the fold point of the layer equation. This leads to  
\begin{align*}
\frac{d\varphi}{dt} & =\omega+\mu-\sin(\varphi^*(\mu))-\cos(\varphi^{*}(\mu))\left(\varphi-\varphi^{*}(\mu)\right),\\
\frac{d\mu}{dt} & =\epsilon\left(-\mu+\eta\left(1-
\sin(\varphi^*(\mu)+\alpha)- \cos(\varphi^{*}(\mu)+\alpha)\right)\left(\varphi-\varphi^{*}(\mu)\right)\right).
\end{align*}
Taking into account that $\varphi^*(\mu)$ is the critical manifold, 
%we take into account only 
the leading terms in both equations read as:
\begin{align*}
\frac{d\varphi}{dt} & =-\cos(\varphi^{*}(\mu))\left(\varphi-\varphi^{*}(\mu)\right),\\
\frac{d\mu}{dt} & =\epsilon\left(-\mu+\eta\left(1-
\sin(\varphi^*(\mu)+\alpha)\right)\right).
\end{align*}

Denoting $\Delta=\varphi_{u}-\varphi_{l}$, we have
\begin{align*}
\frac{d \Delta}{dt} & =-\cos(\varphi^{*}(\mu))\Delta(t),\\
\frac{d \mu}{d t} & =\epsilon\left(-\mu+\eta\left(1-\sin(\varphi^{*}(\mu)+\alpha)\right)\right). 
\end{align*}
The above equations give $d\Delta/dt$ and $d\mu/dt$. Hence, we obtain $d\Delta/d\mu$ as follows
\[
\frac{d\Delta}{d\mu}=\frac{1}{\epsilon}\frac{-\cos(\varphi^{*}(\mu))}{-\mu+\eta\left(1-\sin(\varphi^{*}(\mu)+\alpha)\right)}\Delta,
\]
which can be solved as 
\begin{equation}
\label{eq:SM-scaling}
\Delta(\mu)=
\Delta(\mu_{0})\exp\left[\frac{1}{\epsilon}\int_{\mu_{0}}^{\mu}\frac{\cos(\varphi^{*}(\mu))}{\mu-\eta\left(1-\sin(\varphi^{*}(\mu)+\alpha)\right)}d\mu\right]
=
\Delta(\mu_{0})\exp\left[ - \varepsilon^{-1} C(\mu,\mu_0)\right]
,
\end{equation}
which provides the scaling as in equation \eqref{eq:scaling} of the main paper. 

In addition, according to \cite{KrupaSzmolyanExtending}, there is only an algebraic contraction of the singular funnel stripe across the fold point. Hence, the exponential estimate \eqref{eq:SM-scaling} remains for the further motion along the fast flow. 

\section{Two adaptively coupled phase rotators (3D model)}

Two adaptively coupled phase rotators are described by the system (8-9) form the main paper with $N=2$. For convenience, we repeat it here:
\begin{equation}\label{eq:2_osc}
    \begin{aligned}
        \frac{d\varphi_1}{dt}&=\omega_1 + \mu - \sin\varphi_1 + \frac{\kappa}{2}\sin(\varphi_2-\varphi_1), \\
        \frac{d\varphi_2}{dt}&=\omega_2 + \mu - \sin\varphi_2 + \frac{\kappa}{2}\sin(\varphi_1-\varphi_2),\\
        \frac{d\mu}{dt}&=\varepsilon(-\mu +\eta(1-X)), 
    \end{aligned}
\end{equation}
where $X=\frac{1}{2}\left(\sin(\varphi_1+\alpha) +\sin(\varphi_2+\alpha) \right).$

\begin{figure}
    \centering
       \includegraphics[width=0.49\linewidth]{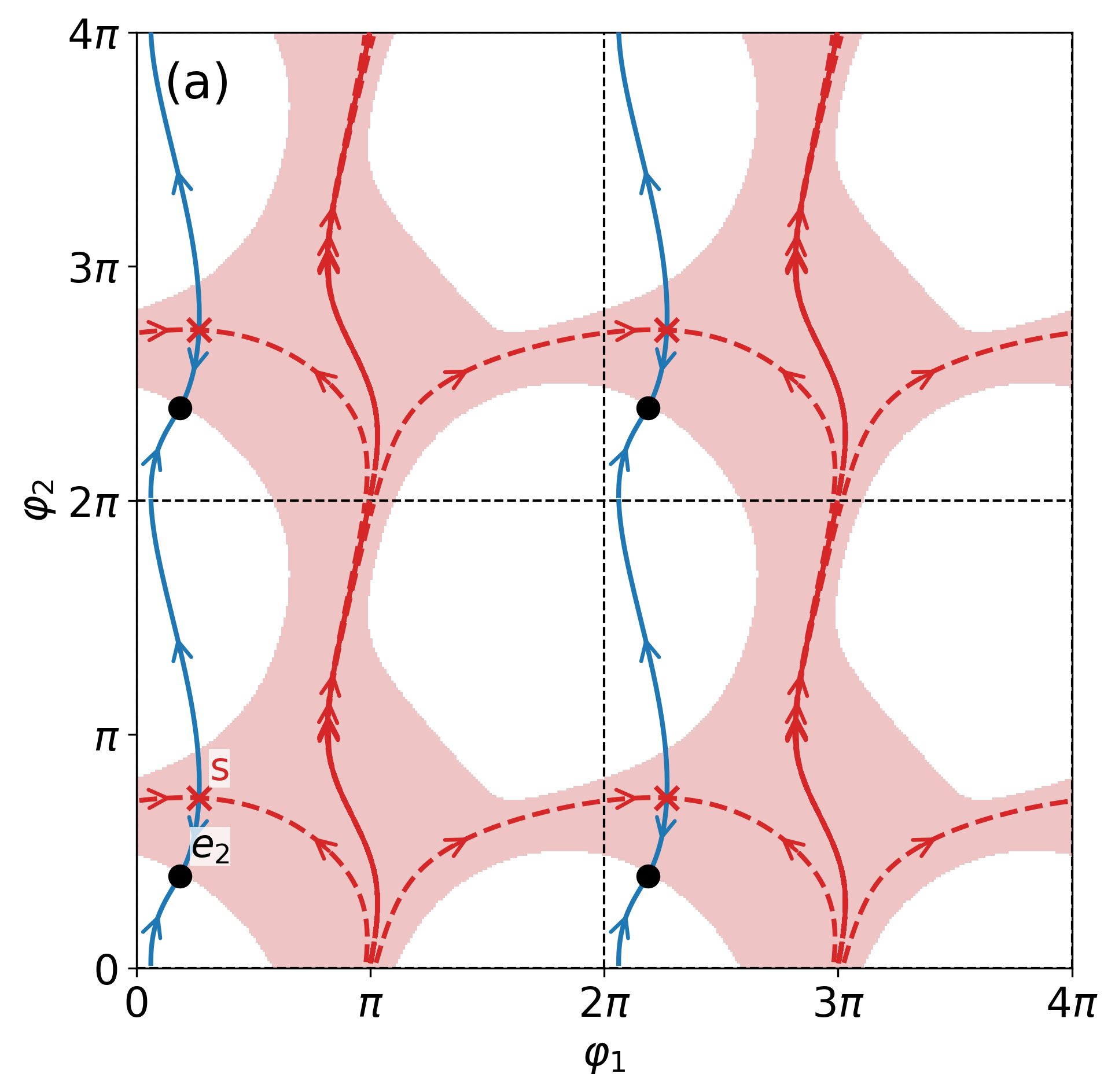}
    \includegraphics[width=0.49\linewidth]{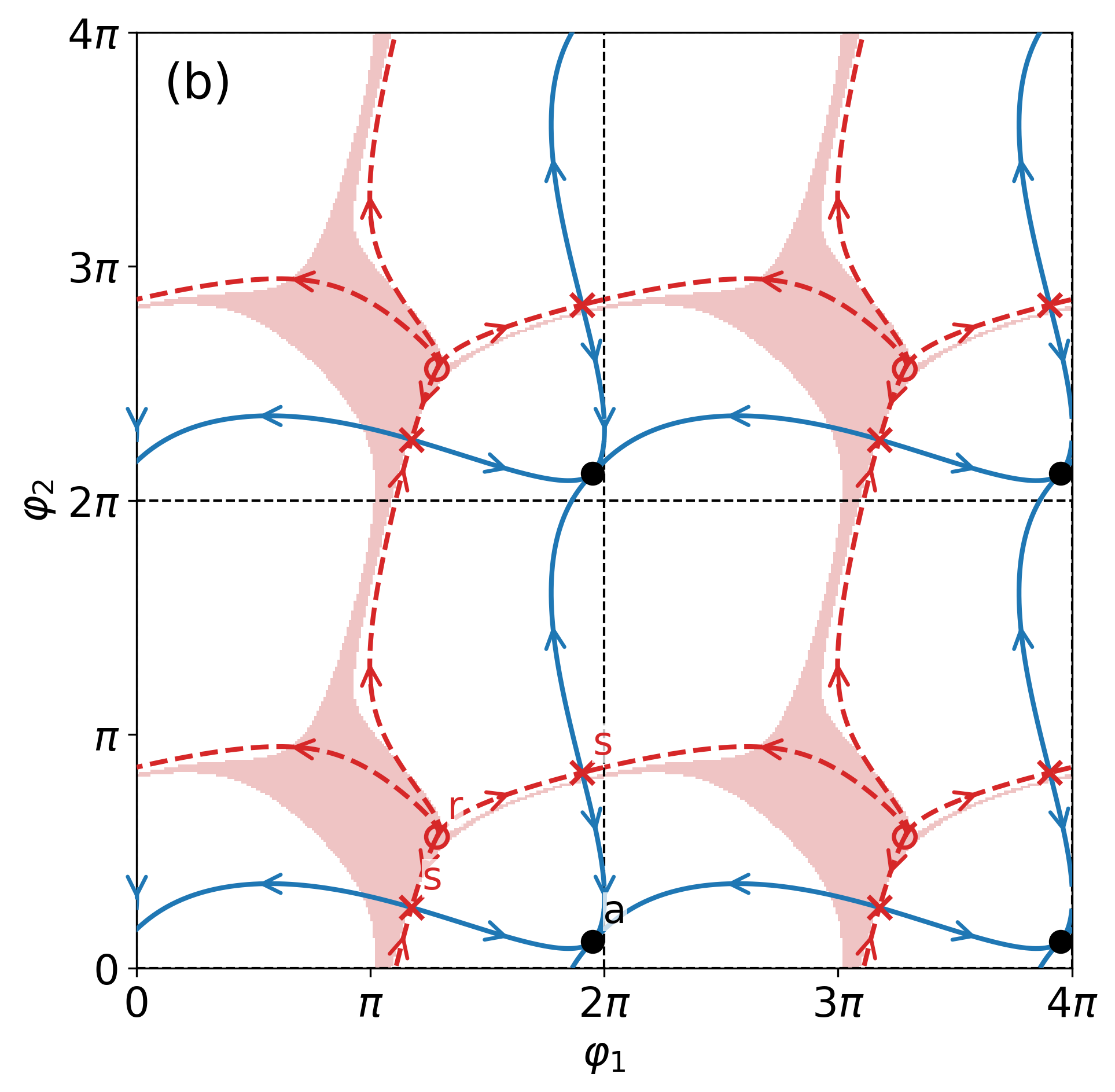}
    \caption{Phase portrait of the layer systems for \eqref{eq:2_osc}, i.e., with a fixed $\mu$, superimposed with the cross-section of the singular basin with the same fixed $\mu$ value. 
    Black dots denote stable equilibria and crosses denote the saddle equilibria of the layer system. The red dashed lines and crosses form an unstable invariant set of the layer system, which gives rise to an unstable critical manifold when $\mu$ is parametrically changed. 
    (a): $\mu=4.25$. This parameter value corresponds to the $\mu$ value for the equilibrium $e_2$ of the full 3D system, hence, this equilibrium also appears in this layer phase portrait. 
    (b) $\mu=3.6$.
    Other parameters: $\omega_1 = -4$, $\omega_2 = -3$, $\kappa = 1$, $\eta = 10$, $\alpha = \pi/2$, $\varepsilon = 0.05$. }
    \label{fig:2acpr}
\end{figure}

The cross-sections of the singular basin for fixed values of $\mu$ are shown in Fig.~\ref{fig:2acpr} as red regions. These cross-sections are superimposed on the phase portraits of the fast subsystem at the corresponding fixed values of $\mu$. It can be seen that the singular basin tends to align with the unstable invariant sets of the fast subsystem (red lines and points).  
In Fig.~\ref{fig:2acpr}(a), for $\mu=4.25$, this unstable invariant set is composed of a repelling rotating limit cycle, a saddle equilibrium $s$ (red crosses), and the stable manifolds of the saddle (dashed red trajectories). 
In Fig.~\ref{fig:2acpr}(b), for $\mu=3.6$, this set is composed of a repelling equilibrium $r$ (red open circles), two saddle equilibria $s$ (red crosses), and the stable manifolds of the saddle (dashed red trajectories).
These unstable invariant sets represent cross-sections of the unstable critical set of the full system \eqref{eq:2_osc}, which plays the role similar to the simple unstable critical manifold $x=0$ of the 2D model \eqref{eq:normal-form-a} from the main manuscript or the unstable part of $S$ from the single rotator model (4-5). Specifically, the basin boundary, i.e., the stable manifold of $e_2$, converges to this set as time goes backwards. More exactly, it converges to the corresponding unstable slow set, which is a perturbation of the critical set.
The structure of this unstable critical set determines the star-shaped form of the cross-section of the basin shown in Fig.~3(d) of the main manuscript.

\makeatletter
\renewcommand{\@biblabel}[1]{[S#1]}
\makeatother

\end{document}